\documentclass[12pt]{article}
\usepackage{array}
\usepackage{ifthen}
\usepackage{url}
\usepackage{float}
\nonstopmode

\usepackage[left=2cm,right=2cm,top=2cm,bottom=2cm]{geometry}
\usepackage[utf8]{inputenc}
\usepackage[english]{babel}
\usepackage[T1]{fontenc}
\usepackage{amsmath}
\usepackage{amsthm}
\usepackage{amsfonts}
\usepackage{amssymb}
\usepackage{color}
\usepackage{multicol}
\usepackage{hyperref}
\hypersetup{
colorlinks=true, 
breaklinks=true, 
urlcolor= blue, 
linkcolor= blue, 
citecolor=blue, 
}
\usepackage{graphicx}
\usepackage{lscape}
\usepackage{setspace}

\newtheorem{theorem}{Theorem}[section]
\newtheorem{lemma}{Lemma}[section]

\def\KK{\mathbb K}
\def\QQ{\mathbb Q}
\def\ZZ{\mathbb Z}

\newcommand{\klammern}[4][]%
{\ifthenelse{\equal{#1}{}}{\left#2}{\csname#1\endcsname#2}%
#4\ifthenelse{\equal{#1}{}}{\right#3}{\csname#1\endcsname#3}}
\newcommand{\betrag}[2][]{\klammern[#1]{\lvert}{\rvert}{#2}}

\newcommand{\conj}[1]{^{(#1)}}

\newcommand{\abs}[1]{\left\vert#1\right\vert}

\renewcommand{\epsilon}{\DONOTUSE}  

\renewcommand{\vartheta}{\DONOTUSE} 
\title{On Pillai's problem involving two linear recurrent sequences : Padovan and Fibonacci} 

\author{ Pagdame Tiebekabe and Serge Adonsou\\}
\begin{document}
\maketitle
\begin{abstract}
In this paper, we find all integers $c$ having at least two representations as a difference between two linear recurrent sequences.
\end{abstract}


\maketitle

\section{Introduction}\label{ppsec1}

It is well-known that the sequence $\{\mathcal{P}_k\}_{k\geq1}$ of Padovan numbers is defined by
$$
\mathcal{P}_0=\mathcal{P}_1=\mathcal{P}_2=1,\quad \mathcal{P}_{k+3}=\mathcal{P}_{k+1}+\mathcal{P}_k,\quad k\geq0.
$$
The first Padovan numbers are 
$$
1, 1, 1, 2, 2, 3, 4, 5, 7, 9, 12, 16, 21, 28, 37, 49, 65, 86, 114, 151, 200, 265  \ldots
$$
The sequence $\{F_k\}_{k\geq1}$ of Fibonacci numbers is defined by
$$
F_0=0,\quad F_1=1,\quad F_{k+2}=F_{k+1}+F_k,\quad k\geq0.
$$
The first Fibonacci numbers are 
$$
0, 1, 1, 2, 3, 5, 8, 13, 21, 34, 55, 89, 144, 233, 377, 610, 987, \ldots
$$
In this paper, we are interested in the Diophantine equation
\begin{equation}\label{ppeq:main}
\mathcal{P}_m-F_n=c
\end{equation}
for a fixed $c$ and variable $m$ and $n$. In particular, we are interested in those integers $c$ admitting at least two representations as a difference between 
a Padovan number and Fibonacci number. This is a variation of the equation 
\begin{equation}\label{ppeq2}
a^x-b^y=c,
\end{equation}
in non-negative integers $(x,y)$ where $a, b, c$ are given fixed positive integers.  The history of equation \eqref{ppeq2} is very rich and goes back to 1935 when Herschfeld \cite{Herschfeld:1935}, \cite{Herschfeld:1936} studied the particular case $(a, b)=(2, 3)$. Extending Herschfeld's work, Pillai \cite{Pillai:1936}, \cite{Pillai:1937} 
proved that if $a, b$ are coprime positive integers then there exists $c_0(a, b)$ such that if $c>c_0(a, b)$ is an integer, then equation \eqref{ppeq2} has at most one positive integer solution $(x, y)$. Since then, variations of equation \eqref{ppeq2} has been intensively studied. Some recent results related to equation \eqref{ppeq:main} are obtained by M. Ddamulira and his collaborators in which they replaced Pillai's numbers by the Fibonacci numbers $F_n$ (see \cite{Ddamulira-Luca-Rakotomalala:2017}), Tribonacci numbers  (see \cite{Bravo-Luca-Yazan:2017}), and $k$-generalized Fibonacci numbers (see \cite{Ddamulira-Gomez-Luca:2017}). 

The aim of this paper is to prove the following result. 
\begin{theorem}\label{ppth:principal}
	The only integers $c$ having at least two representations of the form $\mathcal{P}_m-F_n$ with $m>3, n>1$ are $$\begin{array}{c}
	c \in\{-226, -82, -52, -30, -27, -18, -9, -6, -5, -4, -3, -1, 0, 1, 2, 3, 4, 6,\\ 7, 8, 10, 11, 13, 15, 16, 20, 25, 31, 32, 36, 44, 52, 62, 111, 262\}
	\end{array}.$$ 
\end{theorem}

We organize this paper as follows. In Section~\ref{ppsec2}, we recall some results useful for the proof of Theorem \ref{ppth:principal}. 
The proof of Theorem \ref{ppth:principal} is done in the last section.

\section{Auxiliary results}\label{ppsec2}
\subsection{Some properties of Fibonacci and Padovan sequences }
Here we recall a few properties of  the Fibonacci sequence $\{F_k\}_{k\geq 0}$ and Padovan sequences $\{\mathcal{P}_k\}_{k\geq 0}$  which are useful to proof our theorem.\\
The characteristic equation of Padovan sequence is 
$$
x^3-x-1=0,
$$
has roots $\alpha, \beta, \gamma=\overline{\beta}$, where
$$
\alpha=\dfrac{r_1+r_2}{6},\quad \beta=\dfrac{-r_1-r_2+i\sqrt{3}(r_1-r_2)}{12},
$$
and
$$
r_1=\sqrt[3]{108+12\sqrt{69}}\;\text{ and }\; r_2=\sqrt[3]{108-12\sqrt{69}}.
$$
Further, Binet's formula is
\begin{equation}\label{ppeq:Binet_Padovan}
\mathcal{P}_k=a\alpha^k+b\beta^k+c\gamma^k,\; \text{ for all } k\geq 0,
\end{equation}
where

\begin{equation}\label{ppeq:Formula_a_b_c}
\begin{array}{lll}
a &=& \dfrac{(1-\beta)(1-\gamma)}{(\alpha-\beta)(\alpha-\gamma)}= \dfrac{1+\alpha}{-\alpha^2+3\alpha+1},\vspace{1mm}\\
b &=& \dfrac{(1-\alpha)(1-\gamma)}{(\beta-\alpha)(\beta-\gamma)}=\dfrac{1+\beta}{-\beta^2+3\beta+1},\vspace{1mm}\\
c &=& \dfrac{(1-\alpha)(1-\beta)}{(\gamma-\alpha)(\gamma-\beta)}=\dfrac{1+\gamma}{-\gamma^2+3\gamma+1}=\overline{b}.
\end{array}
\end{equation}
Numerically, we have
\begin{equation}\label{ppeq:Numericle_alpha_beta_a_b}
\begin{array}{l}
1.32<\alpha<1.33,\\
0.86<|\beta|=|\gamma|=\alpha^{-1/2}<0.87,\\
0.72<a<0.73,\\
0.24<|b|=|c|<0.25.
\end{array}
\end{equation}
Using induction, we can prove that
\begin{equation}\label{ppeq:estimate_Padovan}
\alpha^{k-2}\leq \mathcal{P}_k \leq \alpha^{k-1},
\end{equation}
for all $k\geq 4$. 

On the other hand, let $(\delta,
\eta)=\left(\frac{1+\sqrt{5}}{2},\frac{1-\sqrt{5}}{2}\right)$ be the roots of the characteristic equation $x^2-x-1 = 0$ of the Fibonacci sequence $\{F_k\}_{n\geq 0}$. The Binet formula for $F_k$
\begin{equation}\label{ppeq:Binet-Fibonacci}
F_k=\dfrac{\delta^k-\eta^k}{\sqrt{5}} \quad \text{holds for all } k\geq 0.
\end{equation}
This implies easily that the inequality 
\begin{equation}\label{ppeq:enca-Fibonacci}
\delta^{k-2} \leq F_k \leq \delta^{k-1}
\end{equation}
holds for all positive integers $k.$
\subsection{A lower bound for linear forms in logarithms}
The next tools are related to the transcendental approach to solve Diophantine equations. For a non-zero algebraic number $\gamma$ of degree $d$ over $\QQ$, whose minimal polynomial
over $\ZZ$ is $a\prod_{j=1}^d \left(X-\gamma\conj j \right)$, we denote by
\[
h(\gamma) = \frac{1}{d} \left( \log|a| + \sum_{j=1}^d \log\max\left(1,
\betrag{\gamma\conj j}\right)\right)
\]
the usual absolute logarithmic height of $\gamma$.
\begin{lemma}\label{pplem:Matveev}
Let $\gamma_1,\ldots ,\gamma_s$ be a real algebraic numbers and let  $b_1,\ldots,b_s$ be nonzero rational integer numbers. Let $D$ be the degree of the number field $\QQ(\gamma_1,\ldots ,\gamma_s)$ over $\QQ$ and let $A_j$ be a positive real number satisfying  
$$
A_j=\max\{Dh(\gamma_j),|\log\gamma_j|,0.16\}\;\text{for}\; j=1,\ldots ,s. 
$$
Assume that
$$
B\geq\max\{|b_1|,\ldots,|b_s|\}.
$$
If $\gamma_1^{b_1}\cdots\gamma_s^{b_s}\neq1$, then
$$
|\gamma_1^{b_1}\cdots\gamma_s^{b_s}-1|\geq\exp(-C(s,D)(1+\log B) A_1\cdots A_s),
$$
where $C(s,D):=1.4\cdot 30^{s+3}\cdot s^{4.5}\cdot D^2(1+\log D).$
\end{lemma} 

\subsection{A generalized result of Baker-Davenport}
\begin{lemma}\label{pple:Baker-Davenport}
Assume that $\tau$ and $\mu$ are real numbers and $M$ is a positive integer. Let $p/q$ be the convergent of the continued fraction of the irrational $\tau$ such that $q>6M$, and let $A,B,\mu$ be some real numbers with $A>0$ and $B>1$. Let $\varepsilon=||\mu q||-M\cdot||\tau q||$, where $||\cdot||$ denotes the distance from the nearest integer. If    $\varepsilon>0$, then there is no solution of the inequality
$$
0 < m \tau - n + \mu <A B^{-k}
$$
in positive integers $m$, $n$ and $k$ with 
$$
m\leq M\quad \text{ and }\quad k\geq\dfrac{\log(Aq/\varepsilon)}{\log B}.
$$ 
\end{lemma}

\section{Proof of Theorem \ref{ppth:principal}}\label{ppsec3}

Assume that there exist positive integers $n,m,n_1,m_1$ such that $(n,m)\neq (n_1,m_1),$ and 
$$
F_n-\mathcal{P}_m=F_{n_1}-\mathcal{P}_{m_1}.
$$
Because of the symmetry, we can assume that $m \geq m_1$. If $m=m_1$, then $F_n=F_{n_1}$, so $(n, m)=(n_1,m_1)$, contradicting our assumption. Thus, $m > m_1$. Since
\begin{equation}\label{ppeq:Fibonacci_2}
F_n-F\\
_{n_1}=\mathcal{P}_m-\mathcal{P}_{m_1},
\end{equation}
and the right-hand side is positive, we get that the left-hand side is also positive and so $n>n_1$. Thus, $n\geq 2$ and $n_1\geq 1$. Using the Binet's formulas \eqref{ppeq:Binet-Fibonacci} and \eqref{ppeq:Binet_Padovan}, the equation \eqref{ppeq:Fibonacci_2} implies that 
\begin{subequations}
\begin{align}\label{ppeq:maj_n}
{\delta}^{n-4}\leq F_{n-2} \leq F_{n}-F_{n_1}=\mathcal{P}_m-\mathcal{P}_{m_1}<\alpha^{m-1},\\\label{ppeq:min_n}
{\delta}^{n-1}\geq F_{n} > F_{n}-F_{n_1}=\mathcal{P}_m-\mathcal{P}_{m_1}=\mathcal{P}_{m-5}\geq \alpha^{m-7},
\end{align}
\end{subequations}
therefore 
\begin{equation}\label{ppeq:bound_n}
1+\left(\dfrac{\log\alpha}{\log{\delta}}\right) (m-1) < n < \left(\dfrac{\log\alpha}{\log{\delta}}\right) (m-7)+4,
\end{equation}
where $\dfrac{\log\alpha}{\log{\delta}}=0.5843\ldots.$
If $n<300$, then $m\leq 190$. We ran a computer program for $2 \leq n_1 < n \leq 300$ and $1 \leq m_1 < m < 190$ and found only the solutions listed in the \eqref{pplist} at the end of the paper. From now, we assume that $n \geq 300$.\\
Note that the inequality \eqref{ppeq:bound_n} implies that $m <2n$. So, to solve equation \eqref{ppeq:Fibonacci_2}, we need an upper bound for $n$.

\subsection{Bounding $n$}
Note that using the numerical inequalities \eqref{ppeq:Numericle_alpha_beta_a_b} we have 
\begin{equation}\label{ppnumer_ineq}
	\dfrac{|\eta|^n}{\sqrt{5}} + \dfrac{|\eta|^{n_1}}{\sqrt{5}} + |b| |\beta|^m +|c| |\gamma|^m + |b| |\beta|^{m_1} +|c| |\gamma|^{m_1}<1.9.
\end{equation}
Using the Binet formulas in the Diophantine equation \eqref{ppeq:Fibonacci_2}, we get
$$
\begin{array}{rcl}
\abs {\dfrac{{\delta}^n}{\sqrt{5}}-a\alpha^m} & = & \abs {\dfrac{\eta^n}{\sqrt{5}}+\dfrac{{\delta}^{n_1}-\eta^{n_1}}{\sqrt{5}}+(b\beta^m+c\gamma^m) - (a\alpha^{m_1}+b\beta^{m_1}+c\gamma^{m_1})}\vspace{1mm}\\
  &\leq& \dfrac{{\delta}^{n_1}}{\sqrt{5}}+ a \alpha^{m_1} + \dfrac{|\eta|^n}{\sqrt{5}} + \dfrac{|\eta|^{n_1}}{\sqrt{5}} + |b| |\beta|^m +|c| |\gamma|^m + |b| |\beta|^{m_1} +|c| |\gamma|^{m_1} \vspace{1mm}\\
  &<& \dfrac{{\delta}^{n_1}}{\sqrt{5}}+ a \alpha^{m_1}+1.9 \vspace{1mm}\\
   &<& 3.08 \max\{{\delta}^{n_1},\alpha^{m_1}\}.
\end{array}
$$
Dividing through by $a\alpha^m$ and using the relation \eqref{ppeq:maj_n}, we obtain
$$
\begin{array}{lll}
\abs {(\sqrt{5}a)^{-1} {\delta}^n \alpha^{-m}-1}  &<&  \max\left\{ \dfrac{3.08}{a\alpha^m}{\delta}^{n_1},\dfrac{3.08}{a}\alpha^{m_1-m}\right\}\vspace{1mm}\\
  &<&  \max\left\{3.24 \dfrac{{\delta}^{n_1}}{ {\delta}^{n-4}},4.28\alpha^{m_1-m}\right\}.
\end{array}
$$
Hence, we get
\begin{equation}\label{ppeq:Lambda}
\abs {(\sqrt{5}a)^{-1} {\delta}^n \alpha^{-m}-1}   < \max\{{\delta}^{n_1-n+6},\alpha^{m_1-m+3}\}.
\end{equation}
For the left-hand side, we apply Theorem \ref{pplem:Matveev} with the data
$$
s=3, \quad \gamma_1=\sqrt{5}a, \quad \gamma_2=\delta, \quad \gamma_3=\alpha, \quad b_1=-1, \quad b_2=n, \quad b_3=-m.
$$
Throughout we work with $\KK:=\QQ(\sqrt{5},\alpha)$ with $D=6$. Since $\max\{1, n, m\}\leq 2n$ we take $B:=2n$. We have 
$$
h(\gamma_2)=\dfrac{\log{\delta}}{2} \quad \text{and} \quad h(\gamma_3)=\dfrac{\log\alpha}{3}.
$$
Further, the minimal polynomial of $\gamma_1$ is $529x^6-1265x^4-250x^2-125$, then
$$
h(\gamma_1)\approx 1.204.
$$  
Thus, we can take
$$
A_1=7.23, \quad A_2=3\log{\delta}, \quad A_3= 2\log\alpha.
$$
Put
$$
\Lambda=(\sqrt{5}a)^{-1} {\delta}^n \alpha^{-m}-1.
$$
If $\Lambda=0$, then ${\delta}^n(\alpha^{-1})^m=\sqrt{5}a$, which is false, since ${\delta}^n(\alpha^{-1})^m\in\mathcal{O}_{\KK}$ whereas $\sqrt{5}a$ does not, as can be observed immediately from its minimal polynomial. Thus, $\Lambda\neq 0$. Then, by Lemma \ref{pplem:Matveev}, the left-hand side of \eqref{ppeq:Lambda} is bounded as
$$
\log\abs\Lambda > -1.4\cdot 30^6 \cdot 3^{4.5} \cdot 6^2 (1+\log 6) (1+\log 2n) (7.23) (3\log{\delta}) (2\log\alpha).
$$
Comparing with \eqref{ppeq:Lambda}, we get
$$
\min\{ (n-n_1-6) \log {\delta}, (m-m_1-3) \log\alpha \} < 8.45\times 10^{13} (1+\log 2n),
$$
wich gives
$$
\min\{ (n-n_1) \log {\delta}, (m-m_1) \log\alpha \} < 8.45\times 10^{13} (1+\log 2n).
$$
Now the argument splits into two cases.

\medskip

\textbf{Case 1.} $\min\{ (n-n_1) \log {\delta}, (m-m_1) \log\alpha \}= (n-n_1) \log {\delta}$.

\medskip

In this case, we rewrite \eqref{ppeq:Fibonacci_2} as
$$
\abs{\left(\dfrac{{\delta}^{n-n_1}-1}{\sqrt{5}}\right){\delta}^{n_1}-a\alpha^m} = \abs{-a\alpha^{m_1}+\dfrac{\eta^n}{\sqrt{5}}- \dfrac{\eta^{n_1}}{\sqrt{5}} +(b\beta^m+c\gamma^m) - (b\beta^{m_1}+c\gamma^{m_1})}
$$
by using \eqref{ppnumer_ineq} and dividing by $\alpha^m$, we obtain
\begin{equation}\label{eq:Lambda_1}
\abs {\left( \dfrac{{\delta}^{n-n_1}-1}{\sqrt{5}a} \right){\delta}^{n_1}\alpha^{-m} - 1} <  3.65 \alpha^{m_1-m}.
\end{equation}
We put
$$
\Lambda_1=\left( \dfrac{{\delta}^{n-n_1}-1}{a\sqrt{5}} \right){\delta}^{n_1}\alpha^{-m} - 1.
$$
Clearly, $\Lambda_1\neq0$, for if $\Lambda_1=0$, then $ {\delta}^{n}-{\delta}^{n_1}=\sqrt{5}a\alpha^m$. This is impossible if $\sqrt{5}a\alpha^m\in \QQ(\sqrt{5},\alpha)$ but $\not\in\QQ(\sqrt{5})$. Therefore, let us assume that $\sqrt{5}a\alpha^m\in\QQ(\sqrt{5})$. Since $a\alpha^m\in\QQ(\alpha)$ and $\QQ(\alpha)\cap \QQ(\sqrt{5})=\QQ$, we deduce from $\sqrt{5}a\alpha^m\in\QQ(\sqrt{5})$ that we have $\sqrt{5}a\alpha^m=y\sqrt{5}$ for some $y\in\QQ$. Let $\sigma \neq id$ be the unique non trivial $\QQ$-automorphism over $\QQ(\sqrt{5})$. Then, we get
$$
{\delta}^n-{\delta}^{n_1}=\sqrt{5}a\alpha^m=y\sqrt{5}=-\sigma(\sqrt{5}a\alpha^m)=-\sigma({\delta}^n-{\delta}^{n_1})=\eta^{n_1}-\eta^n.
$$
However, the absolute value of the left-hand side is at least ${\delta}^n-{\delta}^{n_1} \geq {\delta}^{n-2} \geq {\delta}^{...}>2$, while the absolute value of right-hand side is at most $\abs{\eta^{n_1}-\eta^n} \leq |{\eta}|^{n_1}+|{\eta}|^n <2.$ By this obvious contradiction we conclude that $\Lambda_1\neq 0$.

We apply Lemma \ref{pplem:Matveev} by taking $s=3$, and
$$
\gamma_1=\dfrac{{\delta}^{n-n_1}-1}{\sqrt{5}a}, \quad \gamma_2=\delta, \quad \gamma_3=\alpha, \quad b_1=1, \quad b_2=n_1, \quad b_3=-m.
$$

On the other hand, the minimal polynomial of $a$ is $23x^3-23x^2+6x-1$ and has roots $a$, $b$, $c$. Since $|b| = |c| < 1$ and $a<1$, then $h(a)=\dfrac{\log 23}{3}$.

Thus, we obtain 
\begin{equation}\label{ppmaj_h(gamma_1)}
\begin{array}{rcl}
h(\gamma_1) & \leq & h\left(\dfrac{{\delta}^{n-n_1}+1}{\sqrt{5}} \right)+ h(a)\vspace{2mm}\\
  &\leq & (n-n_1)h\left({\delta}\right)+h(\sqrt{5})+ h(a)+\log(2)\vspace{2mm}\\
  &<&    \dfrac{1}{2} (n-n_1)\log \delta+\log(\sqrt{5})+\dfrac{\log 23}{3}+\log(2)\vspace{2mm}\\
 &<&     4.22 \times 10^{13}\cdot  (1+\log 2n).
\end{array}
\end{equation}
So, we can take $A_1:=2.53 \times 10^{14} (1+\log 2n)$. Further, as before, we can take $A_2:=3\log{\delta}$ and $A_3:=2\log\alpha$. Finally, since $\max\{1,n_1,m\}\leq 2n $, we can take $B:=2n$.
We then get that
$$
\log\abs{\Lambda_1} > -1.4\cdot 30^6 \cdot 3^{4.5} \cdot 6^2 (1+\log 6) (1+\log 2n) \times (2.53 \times 10^{14} (1+\log 2n)) (3\log{\delta}) (2\log\alpha).
$$
Thus,
$$
\log\abs{\Lambda_1}> -2.96\cdot 10^{27} (1+\log 2n)^2.
$$
Comparing this with \eqref{eq:Lambda_1}, we get that
$$
(m-m_1)\log \alpha<2.96\cdot 10^{27} (1+\log 2n)^2.
$$

\medskip
\textbf{Case 2.} $\min\{ (n-n_1) \log {\delta}, (m-m_1) \log\alpha \}= (m-m_1) \log\alpha $.

\medskip

In this case, we rewrite \eqref{ppeq:Fibonacci_2} as 
$$
\abs{\dfrac{{\delta}^{n}}{\sqrt{5}}-a\alpha^m+a\alpha^{m_1}} = \abs{\dfrac{\eta^n}{\sqrt{5}}+ \dfrac{{\delta}^{n_1}-\eta^{n_1}}{\sqrt{5}} +(b\beta^m+c\gamma^m) - (b\beta^{m_1}+c\gamma^{m_1})}
$$
so 
\begin{equation}\label{ppeq:Lambda_2}
\abs {\dfrac{{\delta}^n\alpha^{-m_1}}{\sqrt{5} a (\alpha^{m-m_1}-1)} - 1} < \dfrac{2.35}{\sqrt{5}a(1-\alpha^{m_1-m})\alpha} \dfrac{{\delta}^{n_1}}{\alpha^{m-1}} < 
17 {\delta}^{n_1-n+4}.
\end{equation}
Let 
$$
\Lambda_2=(\sqrt{5}a(\alpha^{m-m_1}-1))^{-1} {\delta}^{n}\alpha^{-m_1} - 1.
$$
Clearly, $\Lambda_2\neq0$, for if $\Lambda_2=0$ implies ${\delta}^{2n}=5\alpha^{2m_1}a^2(\alpha^{m-m_1}-1)^2.$,However, ${\delta}^{2n}\in\QQ(\sqrt{5})\backslash\QQ$, whereas $5\alpha^{2m_1}a^2(\alpha^{m-m_1}-1)^2\in\QQ(\alpha)$, which is not possible. 

We apply  again Lemma \ref{pplem:Matveev}. In this application, we take again $s=3$, and
$$
\gamma_1=\sqrt{5}a(\alpha^{m-m_1}-1), \quad \gamma_2=\delta, \quad \gamma_3=\alpha, \quad b_1=-1, \quad b_2=n, \quad b_3=-m_1.
$$
We have
$$
\begin{array}{lll}
h(\alpha^{m-m_1}-1)  &\leq& h(\alpha^{m-m_1})+h(-1)+\log 2=(m-m_1)h(\alpha)+\log 2\vspace{1mm}\\
 &=& \dfrac{(m-m_1)\log\alpha}{3} +\log 2 < 9.51 \times 10^{13} (1+\log 2n).
\end{array}
$$
Thus, we obtain
$$
\begin{array}{lll}
h(\gamma_1) &<& 2.82 \times 10^{13} (1+\log 2n)+\dfrac{\log 23}{3}+\log\sqrt{5}\vspace{1mm}\\
  &<& 2.82 \times 10^{13} (1+\log 2n).
\end{array}
$$
So, we can take $A_1:=1.69 \times 10^{14} (1+\log 2n)$. Further, as before, we can take $A_2:=3\log{\delta}$ and $A_3:=2\log\alpha$. Finally, since $\max\{1,n,m_1+1\}\leq 2n $, we can take $B:=2n$.

We then get that
$$
\log\abs{\Lambda_2} > -1.4\cdot 30^6 \cdot 3^{4.5} \cdot 6^2 (1+\log 6) (1+\log 2n) \times (1.69 \times 10^{14} (1+\log 2n)) (3\log{\delta}) (2\log\alpha).
$$
Thus,
$$
\log\abs{\Lambda_1}> -1.97\cdot 10^{27} (1+\log 2n)^2.
$$
Comparing this with \eqref{ppeq:Lambda_2}, we get that
$$
(n-n_1)\log {\delta}<1.97\cdot 10^{27} (1+\log 2n)^2.
$$
Thus, in both Case 1 and Case 2, we have
\begin{subequations}
\begin{align}
\min\{(n-n_1)\log{\delta},(m-m_1)\log\alpha\}<& 8.45\times 10^{13} (1+\log 2n)\label{min}\\
\max\{(n-n_1)\log{\delta},(m-m_1)\log\alpha\}<& 2.96\cdot 10^{27}(1+\log 2n)^2.\label{max}
\end{align}
\end{subequations}

We now finally rewrite equation \eqref{ppeq:Fibonacci_2} as
$$
\abs{\dfrac{{\delta}^{n}}{\sqrt{5}}-\dfrac{{\delta}^{n_1}}{\sqrt{5}}-a\alpha^m+a\alpha^{m_1}} = \abs{\dfrac{\delta^n}{\sqrt{5}}- \dfrac{\delta^{n_1}}{\sqrt{5}} +(b\beta^m+c\gamma^m) - (b\beta^{m_1}+c\gamma^{m_1})} <1.9.
$$
Dividing both sides by $a\alpha^{m_1}(\alpha^{m-m_1}-1)$, we get
\begin{equation}\label{eq:lambda_3}
\abs {\left( \dfrac{{\delta}^{n-n_1}-1}{\sqrt{5}a (\alpha^{m-m_1}-1)}\right) {\delta}^{n_1} \alpha^{-m_1}-1}< \dfrac{5.84}{a(1-\alpha^{m_1-m})\alpha} \dfrac{1}{\alpha^{m-1}}  < 13.8{\delta}^{4-n}. 
\end{equation}
To find a lower-bound on the left-hand side, we use again Lemma \ref{pplem:Matveev} with $s=3$, and 
$$
\gamma_1=\dfrac{{\delta}^{n-n_1}-1}{\sqrt{5}a (\alpha^{m-m_1}-1)}, \quad \gamma_2=\delta,\quad \gamma_3=\alpha,\quad b_1=1,\quad b_2=n_1, \quad b_3= -m_1.
$$
Using $h(x/y)=h(x)+h(y)$ for any two nonzero algebraic numbers $x$ and $y$, we have
$$
\begin{array}{rcl}
h(\gamma_1) & \leq & h\left(\dfrac{{\delta}^{n-n_1}-1}{\sqrt{5}a}\right) +h(\alpha^{m-m_1}-1)\vspace{1mm}\\
  &<& \dfrac{1}{2}(n-n_1) \log {\delta}+\log \sqrt{5} +\dfrac{\log 23}{3} + \dfrac{(m-m_1)\log\alpha}{3} +\log 2\vspace{1mm}\\
            & < & 2.47\cdot 10^{27} (1+\log 2n)^2,
\end{array}
$$
where in the above chain of inequalities, we used the argument from 
\eqref{ppmaj_h(gamma_1)} as well as the bound \eqref{max}. So, we can take $A_1:=1.78\cdot 10^{28} (1+\log 2n)^2$ and certainly $A_2:=3\log{\delta}$ and $A_3:=2\log\alpha$. Using similar arguments as in the proof that $\Lambda_1\neq 0$ we show that if we put
$$
\Lambda_3=\left( \dfrac{{\delta}^{n-n_1}-1}{\sqrt{5}a (\alpha^{m-m_1}-1)}\right) {\delta}^{n_1} \alpha^{-m_1}-1,
$$
then $\Lambda_3\neq 0$. Lemma \ref{pplem:Matveev} gives
$$
\log\abs{\Lambda_3} > -1.4\cdot 30^6 \cdot 3^{4.5} \cdot 6^2 (1+\log 6) (1+\log 2n) \times (1.78\cdot 10^{28} (1+\log 2n)^2) (3\log{\delta}) (2\log\alpha),
$$
which together with \eqref{eq:lambda_3} gives 
$$
(n-4)<2.08\cdot 10^{41} (1+\log 2n)^3,
$$ 
leading to $n<2.83\cdot 10^{47}$.

\subsection{Reducing $n$}
We now need to reduce the above bound for $n$ and to do so we make use of Lemma \ref{pple:Baker-Davenport} several times and each time $M:=2.83\cdot 10^{47}$. To begin with, we return to \eqref{ppeq:Lambda} and put
$$
\Gamma:=n \log{\delta} -m \log \alpha-\log(\sqrt{5}a).
$$
For technical reasons  we assume that $\min\{n-n_1,m-m_1\}\geq 20$. We go back to the inequalities for $\Lambda,~\Lambda_1,~\Lambda_2$.
Since we assume that $\min\{n-n_1,m-m_1\}\geq 20$ we get 
$|e^{\Gamma}-1|=|\Lambda|<\dfrac{1}{4}$. Hence, $|\Lambda|<\dfrac{1}{2}$ and since the inequality $|x|<2|e^x-1|$ holds for all $x\in\left(-\frac{1}{2}, \frac{1}{2}\right)$, we get
$$
\abs\Gamma < 2\max\{{\delta}^{n_1-n+6},\alpha^{m_1-m+3}\}\leq \max\{{\delta}^{n_1-n+8},\alpha^{m_1-m+6}\}.
$$
Assume $\Gamma>0$. We then have the inequality
$$
\begin{array}{rcl}
0<n\left(\dfrac{\log{\delta}}{\log\alpha}\right)-m-\dfrac{\log(1/(\sqrt{5}a))}{\log\alpha} &<& \max\left\{ \dfrac{{\delta}^{8}}{\log\alpha} {\delta}^{-(n-n_1)}, \dfrac{\alpha^{6}}{\log\alpha} \alpha^{-(m-m_1)} \right\} \vspace{2mm}\\
   &<& \max\{ 170\cdot{\delta}^{-(n-n_1)},20\cdot \alpha^{-(m-m_1)} \}.
\end{array}
$$
We apply Lemma \ref{pple:Baker-Davenport} with
$$
\tau=\dfrac{\log{\delta}}{\log{\alpha}},\quad \mu=\dfrac{\log(1/(\sqrt{5}a))}{\log\alpha},\quad (A,B)=(170,\delta)\; \text{ or }\; (20,\alpha).
$$
Let $\tau=[a_0,a_1,\ldots]=[1; 1, 2, 2, 6, 2, 1, 2, 1, 2, 1, 1, 11, 1, 2, 3, 1, 7, 37, 4,\ldots]$ be the continued fraction of $\tau$. We choose consider the 98-th convergent 
$$
\dfrac{p}{q}=\dfrac{p_{98}}{q_{98}} = \dfrac{78093067704223831799032754534503501859635391435517}{45634243076387457097046528084208490147594968308975}.
$$
If satisfied $q=q_{98}>6M$. Further, it yields $\varepsilon>0.35$, and therefore either
$$
n-n_1\leq \dfrac{\log(170q/\varepsilon)}{\log{\delta}}<250,\;\text{ or }\;  m-m_1\leq \dfrac{\log(20q/\varepsilon)}{\log\alpha}< 420.
$$
In the case of $\Gamma<0$, we consider the following inequality instead:
$$
\begin{array}{rcl}
m\left(\dfrac{\log\alpha}{\log{\delta}}\right)-n+\dfrac{\log(\sqrt{5}a)}{\log{\delta}} &<& \max\left\{ \dfrac{{\delta}^{9}}{\log{\delta}} \alpha^{-(n-n_1)}, \dfrac{\alpha^{12}}{\log{\delta}} \alpha^{-(m-m_1)} \right\} \vspace{2mm}\\
   &<& \max\{ 98\cdot {\delta}^{-(n-n_1)},12\cdot \alpha^{-(m-m_1)} \},
\end{array}
$$
instead and apply Lemma \ref{pple:Baker-Davenport} with
$$
\tau=\dfrac{\log\alpha}{\log{\delta}},\quad \mu=\dfrac{\log(\sqrt{5}a)}{\log{\delta}},\quad (A,B)=(98,\delta)\; \text{ or }\; (12,\alpha).
$$
Let $\tau=[a_0,a_1,\ldots]=[0; 1, 1, 2, 2, 6, 2, 1, 2, 1, 2, 1, 1, 11, 1, 2, 3, 1, 7, 37,\ldots]$ be the continued fraction of $\tau$ (note that the current $\tau$ is just the reciprocal of the previous $\tau$). We consider the $98$-th convergent 
$$
\dfrac{p}{q}=\dfrac{p_{98}}{q_{98}} = \dfrac{1000540334879242934726141761162813294034885977722}{1712206861451396832387596141129961335575127483549}.
$$
which satisfies $q=q_{98}>6M$. This yields again $\varepsilon>0.47$, and therefore either
$$
n-n_1\leq \dfrac{\log(98q/\varepsilon)}{\log{\delta}}<242,\;\text{ or }\;  m-m_1\leq \dfrac{\log(12q/\varepsilon)}{\log\alpha}< 406.
$$
These bounds agree with the bounds obtained in the case that $\Gamma>0$. As a conclusion, we have either $n-n_1\leq 250$ or $m-m_1\leq 420$ whenever $\Gamma\neq0$.

Now, we have to distinguish between the cases $n-n_1\leq 250$ and $m-m_1\leq 420$. First, let assume that $n-n_1\leq 250$. In this case, we consider inequality \eqref{eq:Lambda_1} and assume that $m-m_1\geq 20$. We put
$$
\Gamma_1=n_1\log{\delta}-m\log\alpha+\log\left(\dfrac{{\delta}^{n-n_1}-1}{\sqrt{5}a}\right).
$$
Then inequality \eqref{eq:Lambda_1} implies that
$$
\abs{\Gamma_1} <7.3\alpha^{m_1-m}.
$$
If we further assume that $\Gamma_1>0$, we then get
$$
0<n_1\left(\dfrac{\log{\delta}}{\log\alpha}\right)-m+\dfrac{\log(({\delta}^{n-n_1}-1)/(\sqrt{5}a)}{\log\alpha}<26\cdot\alpha^{-(m-m_1)}.
$$
Again we apply Lemma \ref{pple:Baker-Davenport} with the same $\tau$ as in the case when $\Gamma>0$. We use the 98-th convergent ${p}/{q}={p_{98}}/{q_{98}}$ of $\tau$ as before. But in this case we choose $(A,B):=(26,\alpha)$ and use
$$
\mu_k=\dfrac{\log(({\delta}^k-1)/(\sqrt{5}a))}{\log\alpha},
$$
instead of $\mu$ for each possible value of $k:=n-n_1\in [1,2,\ldots 250].$ For the remaining values of $k$, we get
 $\varepsilon>0.0004$. Hence, by Lemma \ref{pple:Baker-Davenport}, we get
$$
m-m_1< \dfrac{\log(26q/0.0004)}{\log\alpha}< 446.
$$ 
Thus, $n-n_1\leq 250$ implies $m-m_1\leq 446$. 

In the case that $\Gamma_1 < 0$ we follow the ideas from the case that $\Gamma_1 > 0$. We use the same $\tau$ as in the case that $\Gamma < 0$ but instead of $\mu$ we take
$$
\mu_k=\dfrac{\log((\sqrt{5}a)/({\delta}^k-1))}{\log{\delta}},
$$
for each possible value of $n- n_1 = k = 1, 2,\ldots , 250$. Using Lemma \ref{pple:Baker-Davenport} with this setting we also obtain in this case that $n - n_1 \leq 250$ implies $m- m_1 \leq 429 $.

Now let us turn to the case that $m-m_1\leq 420$ and let us consider inequality \eqref{ppeq:Lambda_2}. We put
$$
\Gamma_2=n\log{\delta}-m_1\log\alpha+\log(1/(\sqrt{5}a(\alpha^{m-m_1}-1))),
$$
and we assume that $n-n_1 \geq 20$. We then have
$$
\abs{\Gamma_2} < \dfrac{34{\delta}^4}{\alpha^{n-n_1}}.
$$
Assuming $\Gamma_2>0$, we get 
$$
0<n\left(\dfrac{\log{\delta}}{\log\alpha}\right)-m_1+\dfrac{\log((1/(\sqrt{5}a(\alpha^{m-m_1}-1)))}{\log\alpha}< \dfrac{34{\delta}^4}{(\log \alpha)\alpha^{n-n_1}} <830{\delta}^{-(n-n_1)}.
$$
We apply again Lemma \ref{pple:Baker-Davenport} with the same $\tau$, $q$, $M$, $(A,B):=(830,\delta)$ and 
$$
\mu_k=\dfrac{\log((1/(\sqrt{5}a(\alpha^{k}-1)))}{\log \alpha}\quad \text{for }k=1, 2, \ldots 420.
$$
We get $\varepsilon>0.00077$, therefore
$$
n-n_1<\dfrac{\log(830q_{98}/0.00077)}{\log{\delta}}<263.
$$
A similar conclusion is reached when $\Gamma_2<0$. To conclude, we first get that either $n-n_1\leq 250$ or $m-m_1\leq 446$. If $n-n_1\leq 250$, then $m-m_1\leq 446$, and if $m-m_1\leq 420$ then $n-n_1\leq 263$. In conclusion, we always have $n-n_1<263$ and $m-m_1<446$.

Finally we go to \eqref{eq:lambda_3}. We put 
$$
\Gamma_3=n_1\log{\delta}-m_1\log\alpha+\log\left(\dfrac{{\delta}^{n-n_1}-1}{\sqrt{5}a(\alpha^{m-m_1}-1)}\right).
$$
Since $n\geq 300$, inequality \eqref{eq:lambda_3} implies that
$$
\abs{\Gamma_3}<\dfrac{17}{{\delta}^{n-4}}.
$$
Assume that $\Gamma_3>0$. Then
$$
0<n_1\left(\dfrac{\log{\delta}}{\log\alpha}\right)-m_1+\dfrac{\log(({\delta}^k-1)/(\sqrt{5}a(\alpha^l-1)))}{\log\alpha}<390{\delta}^n,
$$
where $(k,l):=(n-n_1,m-m_1)$. We apply again Lemma \ref{pple:Baker-Davenport} with the same $\tau$, $M$, $q$, $(A,B):=(390,\delta)$ and 
$$
\mu_{k,l}=\dfrac{\log(({\delta}^k-1)/(\sqrt{5}a(\alpha^l-1)))}{\log\alpha}\quad \text{for } 1\leq k\leq 264,\; 1\leq l\leq 446.
$$
We consider the $99$th convergent $\dfrac{p_{99}}{q_{99}}$. For all pairs $(k,l)$ we get that $\varepsilon>2\times 10^{-5}$. Thus, Lemma \ref{pple:Baker-Davenport} yields that
$$
n<\dfrac{\log (390 \times q_{99}/\varepsilon)}{\log{\delta}}<274.
$$

Theorem \ref{ppth:principal} is therefore proved.

On the next page is presented the table that gives the couples for which we obtain the different representations of $c$ on the form $\mathcal{P}_m-F_n=c$.

\begin{table}[H]\label{pplist}
	\centering
	\begin{tabular}{|c|c|}
		\hline
		$c$ & $(m,n)$ \\
		\hline
		$-226$& $(8, 13), (19, 14)$  \\
		\hline
		$-82$& $(8, 11), (19, 13)$  \\
		\hline
		$-52$& $(5, 10), (14, 11)$  \\
		\hline
		$-30$& $(6, 9), (18, 12)$  \\
		\hline
		$-27$& $(8, 9), (13, 10)$  \\
		\hline
		$-18$& $(5, 8), (11, 9), (14, 10)$  \\
		\hline
		$-9$& $(6, 7), (10, 8)$  \\
		\hline
		$-6$& $(4, 6), (8, 7), (13, 9), (15, 10)$ \\
		\hline
		$-5$& $(5, 6), (11, 8)$  \\
		\hline
		$-4$& $(6, 6), (9, 7)$  \\
		\hline
		$-3$& $(4, 5), (7, 6), (17, 11)$  \\
		\hline
		$-1$& $(4, 4), (6, 5), (8, 6), (10, 7)$ \\
		\hline
		$0$& $(4, 3), (5, 4), (7, 5), (12, 8)$ \\
		\hline
		$1$& $(4, 2), (5, 3), (6, 4), (9, 6)$ \\
		\hline
		$2$& $(5, 2), (6, 3), (7, 4), (8, 5)$ \\
		\hline
		$3$& $(6, 2), (7, 3), (11, 7), (14, 9)$ \\
		\hline
		$4$& $(7, 2), (8, 4), (9, 5), (10, 6)$ \\
		\hline
		$6$& $(8, 2), (9, 4), (24, 15)$ \\
		\hline
		$7$& $(9, 3), (10, 5), (13, 8), (19, 12)$ \\
		\hline
		$8$& $(9, 2), (11, 6), (12, 7)$ \\
		\hline
		$10$& $(10, 3), (16, 10)$ \\
		\hline
		$11$& $(10, 2), (11, 5)$ \\
		\hline
		$13$& $(11, 4), (12, 6)$ \\
		\hline
		$15$& $(11, 2), (13, 7), (15, 9)$  \\
		\hline
		 16& $(12, 5), (14, 8) $ \\
		\hline 
		20& $(12, 2), (13, 6) $ \\
		\hline
		25& $(13, 4), (18, 11) $ \\
		\hline
		31& $(16, 9), (17, 10) $ \\
		\hline
		32& $(14, 5), (21, 13) $ \\
		\hline
		36& $(14, 2), (15, 7) $ \\
		\hline
		44& $(15, 5), (16, 8) $ \\
		\hline
		52& $(16, 7), (17, 9) $ \\
		\hline
		62& $(16, 4), (19, 11) $ \\
		\hline
		111& $(18, 4), (20, 11) $ \\
		\hline
		262& $(21, 4), (22, 11)$\\
		\hline
	\end{tabular}
	\caption{Representations}
\end{table}

\end{document}